\documentclass[12pt,a4paper]{article}
\usepackage{amssymb}

\usepackage{amsmath}

\newtheorem{theorem}{Theorem}[section]
\newtheorem{definition}[theorem]{Definition}
\newtheorem{lemma}[theorem]{Lemma}
\newtheorem{corollary}[theorem]{Corollary}

\begin{document}
\title{Composition-Diamond Lemma for Modules\footnote{Supported by the
NNSF of China (No.10771077) and the NSF of Guangdong Province
(No.06025062).}}
\author{Yuqun Chen, Yongshan Chen and Chanyan Zhong   \\
%EndAName
\\
{\small \ School of Mathematical Sciences}\\
{\small \ South China Normal University}\\
{\small \ Guangzhou 510631}\\
{\small \ P. R. China}\\
{\small \ yqchen@scnu.edu.cn}\\
{\small \ jackalshan@126.com} \\
{\small \ zhongchanyan@tom.com}}
\date{}
\maketitle

\maketitle \noindent\textbf{Abstract:} In this paper we give some
relationships among the Gr\"{o}bner-Shirshov bases  in free
associative algebras, free left modules and ``double-free" left
modules (free modules over a free algebra). We give the Chibrikov's
Composition-Diamond lemma for modules and show that Kang-Lee's
Composition-Diamond lemma follows from this lemma. As applications,
we also deal with highest weight module over the Lie algebra $sl_2$,
Verma module over a Kac-Moody algebra, Verma module over Lie algebra
of coefficients of a free conformal algebra and the universal
enveloping module for a Sabinin algebra.


\noindent {\bf AMS} Mathematics Subject Classification(2000): 16S15,
13P10, 17A01, 17B67, 16D10

\section{Introduction}

Composition-Diamond lemma for modules was first formulated and
proved by S.-J. Kang and K.-H. Lee in \cite{kl1} and \cite{kl3}.
According to their approach, a Gr\"{o}bner-Shirshov basis of a
cyclic module $M$ over an algebra $A$ is a pair $(S,T)$, where $S$
is the set of the defining relations of $A, \ A=k\langle X|S\rangle$
and $T$ is the defining relations for the $A$-module $_A M=
mod_A\!\langle e|T\rangle$. Then Kang-Lee's Lemma says that $(S,T)$
is a Gr\"{o}bner-Shirshov pair for $A$-module $_A M=mod_A \langle
e|T\rangle$ if $S$ is a Gr\"{o}bner-Shirshov basis of $A$ and $T$ is
closed under the right-justified composition with respect to $S$,
and for $f\in S, \ g\in T$, such that $(f,g)_w$ is defined and
$(f,g)_w\equiv 0 \ mod(S,T;w)$.

They gave applications of this lemma for irreducible modules over
$sl_n(k)$ \cite{kl2}, Specht modules over Hecke algebras and
Ariki-Koike algebras in \cite{kl5} and \cite{kl4}. Some years later,
E. S. Chibrikov \cite{ch} suggested a new Composition-Diamond lemma
for modules that treat any module as a factor module of
``double-free" module, a free module $mod_{k\langle X\rangle}\langle
Y\rangle$ over a free algebra $k\langle X\rangle$. In this approach
any $A$-module $_A M$ is presented in the form
$$
_A M=mod_{k\langle X\rangle} \langle Y|SX^*Y, \ T\rangle,
$$
where $A=k\langle X|S\rangle, \ _A M=mod_A \langle Y|T\rangle,\ X^*$
is the free monoid generated by $X$.

The aim of this paper is to give some relationships among the
Gr\"{o}bner-Shirshov bases  in free associative algebras, free left
modules and ``double-free" left modules.  Also we give some
applications of Composition-Diamond lemma for ``double-free"
modules.

The paper is organized as follows. In \S 2, we are dealing with
Gr\"{o}bner-Shirshov bases and Composition-Diamond lemma for left
ideals of a free algebra. Actually, it is a special case of cyclic
``double-free" modules. In \S 3, we give some relationships among
the Gr\"{o}bner-Shirshov bases  in free associative algebras, free
left modules and ``double-free" modules. We give a proof of
Chibrikov's Composition-Diamond lemma and formulate Kang-Lee's
Composition-Diamond lemma. Then, we show that the latter follows
from the former. As applications, in \S 4, \S 5, \S 6 and \S 7, we
are dealing with highest weight module over the Lie algebra $sl_2$,
Verma module over a Kac-Moody algebra, Verma module over
 Lie algebra of coefficients of a free conformal algebra and the
universal enveloping module for a Sabinin algebra, respectively.

Let $k$ be a field, $X$ a set, $X^*$  the free monoid of associative
words on $X$, and $k\langle X\rangle$ the free associative algebra
over $X$ and $k$. For a word $w\in X^*$, we denote the length of $w$
by $deg(w)$. Suppose that $``<"$ is a well order on $X^*$. For any
polynomial $f$, let $\overline{f}$ be the leading term of $f$. If
the coefficient of $\overline{f}$ is $1$, then we call this
polynomial is monic.

The following lemma will be used in  \S 4, \S 5 and  \S 6.

\begin{lemma}\label{l0}(\cite{s}, cf.  \cite{bm99}, \cite{bkl99}, \cite{lr},  \cite{bc1})
Let $Lie(X)$ be a free Lie algebra over a set $X$ and a field $k$.
Let $S\subset{Lie(X)}$ be a nonempty set of monic Lie polynomials.
Then, with a deg-lex order on $X^*$, $S$ is a Gr\"{o}bner-Shirshov
basis in $Lie(X)$ if and only if $S^{(-)}$ is a Gr\"{o}bner-Shirshov
basis in $k\langle X\rangle$ where $S^{(-)}$ is just $S$ but
substitute all $[xy]$ by $xy-yx$.
\end{lemma}

\section{Composition-Diamond lemma for left ideals of a free algebra}

Let $X$ be a set and $``<"$  a well order on $X^*$. Let $S\subset
k\langle X\rangle$ with each $s\in S$ monic. Then $k\langle X\rangle
S$ is the left ideal of $k\langle X\rangle$ generated by $S$. For
left ideal  $k\langle X\rangle S$, we define the compositions in $S$
as follows.

\begin{definition}
For any  $f, \ g\in S$, if $w=\overline{f}=a\overline{g}$ for some
$a\in X^*$, then the composition of $f$ and $g$ is defined to be
$(f,g)_w=f-ag$. Transformation $f\rightarrow f-ag$ is called the
elimination of the leading word (ELW) of $g$ in $f$.

If $(f,g)_w=\sum \alpha_i a_i s_i$, where $\alpha_i\in k, \ a_i \in
X^*, \ s_i \in S$ and $a_i \overline{s}_i<w$, then we say this
composition is trivial $mod(S,w)$ and denote by  $(f,g)_w\equiv 0 \
\  mod(S,w)$.
\end{definition}

\begin{definition}
Let $S\subset k\langle X\rangle$ with each $s\in S$ monic. Then we
call $S$ a Gr\"{o}bner-Shirshov basis of left ideal $k\langle
X\rangle S$ if all the compositions are trivial modulo $S$. $S$ is
called the minimal Gr\"{o}bner-Shirshov basis of $k\langle X\rangle
S$, if there are no compositions of polynomials in $S$, i.e.
$\overline{f}\neq a\overline{g}$ for any $f, \ g\in S, \ f\neq g$.
\end{definition}

A well order $``<"$ on $X^*$ is left compatible if for any $u, v\in
X^*$,
$$
u > v \Rightarrow wu > wv,  \ for \  all \
 w\in  X^*.
$$
That $``<"$ is right compatible is similarly defined. $``<"$ is
monomial if it is both left and right compatible.

Now we can formulate the Composition-Diamond lemma for left ideals
of a free algebra.

\begin{lemma}(Composition-Diamond lemma for left ideals of $k\langle X\rangle$)\label{l1}
Let $S\subset k\langle X\rangle$ with each $s\in S$ monic and $``<"$
a  left compatible order on $X^*$. Then the following statements are
equivalent:
\begin{enumerate}
\item[(1)]\ $S$ is a Gr\"obner-Shirshov basis of left ideal $k\langle
X\rangle S$.
\item[(2)]\ If \ $0\neq
f\in k\langle X\rangle S$, then $\overline{f}=a\overline{s}$ for
some $a\in X^*, \ s\in S$.
\item[(2')]\  If \ $0\neq f\in k\langle
X\rangle S$, then $f=\sum \alpha_i a_i s_i$ with $a_1
\overline{s}_1>a_2 \overline{s}_2>\cdots$, where each $a_i\in X^*, \
s_i\in S$.
\item[(3)]\
$Red(S)=\{w \in X^*|w\neq a\overline{s}, \ a\in X^*, \ s\in S\}$ is
a $k$-linear basis for the factor $k\langle X\rangle$-module
$_{k\langle X\rangle} k\langle X\rangle/k\langle X\rangle S$.
\end{enumerate}
\end{lemma}

Lemma \ref{l1} is a special case of Lemma \ref{l2} (see the next
section).

Assume that $S$ is a Gr\"obner-Shirshov basis for the left ideal
$k\langle X\rangle S$ with no compositions at all between different
elements of $S$. It means that $S$ is a minimal Gr\"obner-Shirshov
basis of the left ideal $k\langle X\rangle S$. Then, $k\langle
X\rangle S$ is a free $k\langle X\rangle$-module with the basis $S$.

Now, we cite Kang-Lee's Composition-Diamond lemma.

Let $S, \ T\subset k\langle X\rangle, \ A=k\langle X|S\rangle$,
$_AM=_A\!A/A(T+Id(S))$  a left $A$-module and $f, g\in k\langle
X\rangle$. In Kang-Lee's paper \cite{kl1}, the composition of $f$
and $g$ is defined as follows.

\begin{definition}\label{d1} (\cite{kl1},\cite{kl3}) \
\begin{enumerate}
\item[(a)]\ If there exist $a, \ b\in X^*$ such that $w=\overline{f}a=b\overline{g}$,
with $deg(\overline{f})>deg(b)$, then the composition of
intersection is defined to be $(f,g)_w=fa-bg$.
\item[(b)]\ If there exist $a, \ b\in X^*$ such that $w=a\overline{f}b=\overline{g}$,
then the composition of inclusion is defined to be $(f,g)_w=afb-g$.
\item[(c)]\ A composition  $(f,g)_w$ is called right-justified if
$w=\overline{f}=a\overline{g}$ for some $a\in X^*$.
\end{enumerate}
\end{definition}

If $f-g=\sum \alpha_i a_is_i b_i+\sum \beta_j c_j t_j$, where
$\alpha_i, \ \beta_j \in k, \ a_i, \ b_i, \ c_j\in X^*, \ s_i\in S,
\ t_j\in T$ with $a_i \overline{s}_i b_i<w$ and $c_j
\overline{t}_j<w$ for each $i$ and $j$, then $f-g$ is called trivial
with respect to $S$ and $T$ and denote it by $f\equiv g \
mod(S,T;w)$.

When $T=\varnothing$, we simply write $f\equiv g \ mod(S,w)$. If for
any $f,\ g\in S,\ (f,g)_w$ is defined and $f\equiv g \ mod(S,w)$,
then we say $S$ is closed under the composition. Note that, if this
is the case, $S$ is called a Gr\"{o}bner-Shirshov basis in $k\langle
X\rangle$ which is firstly introduced by Shirshov \cite{s} (see also
\cite{b72}, \cite{b76}).

\ \

\noindent{\bf Remark.} If a subset $S$ of $k\langle X\rangle$ is not
a Gr\"{o}bner-Shirshov basis, then we can add to $S$ all nontrivial
compositions of polynomials of $S$, and by continuing this process
(maybe infinitely) many times, we eventually obtain a
Gr\"{o}bner-Shirshov basis $S^{c}$ in $k\langle X\rangle$. Such a
process is called the Shirshov algorithm.

\begin{definition}\label{d4} (\cite{kl1},\cite{kl3}) \ A pair $(S,T)$ of subsets of
monic elements of $k\langle X\rangle$ is called a Gr\"obner-Shirshov
pair if $S$ is closed under the composition, $T$ is closed under the
right-justified composition with respect to $S$ and $T$, and for any
$f\in S, \ g\in T$ and $w\in X^*$ such that if $(f,g)_w$ is defined
(it means that $a\overline{f}b=c\overline{g}$, where $a,b ,c\in X^*,
\ f\in S, \ g\in T$ and $deg(\overline{f})>deg(c)$), we have
$(f,g)_w\equiv 0\ mod \ (S,T;w)$. In this case, we say that $(S,T)$
is a Gr\"obner-Shirshov pair for the $A$-module $_A M=_A
A/A(T+Id(S))$, where $A=k\langle X|S\rangle$.
\end{definition}

The following is the Kang-Lee's Composition-Diamond lemma for a left
module.

\begin{theorem}\label{t1}(\cite{kl1},\cite{kl3})
Let $(S,T)$ be a pair of subsets of monic elements in $k\langle
X\rangle$, let $A=k\langle X|S\rangle$ be the associative algebra
defined by $S$, and let $_AM=_A\!A/A(T+Id(S))$ be the left module
defined by $(S,T)$. Suppose that $(S,T)$ is a Gr\"obner-Shirshov
pair for the $A$-module $_A M$ and $p\in k\langle X\rangle T+Id(S)$.
Then $\overline{p}=a\overline{s}b$ or $\overline{p}=c\overline{t}$,
where $a,b,c\in X^*,  \ s\in S, \ t\in T$.
\end{theorem}

In particular, Lemma \ref{l1} is a special case of Theorem \ref{t1}
 when $S=\varnothing$.

\section{Composition-Diamond lemma for ``double-free"  modules}

Let $X,Y$ be sets and $mod_{k\langle X\rangle} \langle Y\rangle$ a
free left $k\langle X\rangle$-module with the basis $Y$. Then
$mod_{k\langle X\rangle} \langle Y\rangle=\oplus_{y\in Y} k\langle
X\rangle y$ is called a ``double-free"  module.

We now define  Gr\"{o}bner-Shirshov basis in $mod_{k\langle
X\rangle} \langle Y\rangle$.

Suppose that $<$ is a monomial order on $X^*$, $<$ a well order on
$Y$  and $X^*Y=\{uy|u\in X^*, \ y\in Y\}$. We define an order
$``\prec"$ on $X^*Y$: for any $w_1=u_1y_1,w_2=u_2y_2\in X^*Y$,
$$
w_1\prec w_2\Leftrightarrow u_1<u_2 \ \ \mbox{ or } u_1=u_2, y_1<y_2
\ \ \ \  \ \  \  \ \ \  \ \  \  \ \ \ \  \ \  \  \ \ \  \ \  \  \ \
\ (*)
$$
It is clear that the order $``\prec "$ is left compatible in the
sense of
$$
w\prec w'\Rightarrow aw \prec aw' \ \mbox{ for any  } a\in X^*.
$$

 Let $S\subset
mod_{k\langle X\rangle} \langle Y\rangle$ with each $s\in S$ monic.
We define the composition in $S$ only the inclusion composition
which means $w=\bar{f}=a\bar{g}$, where $f,g\in S$.

If $(f,g)_w=f-ag=\sum \alpha_i a_i s_i$, where $\alpha_i\in k, \ a_i
\in X^*, \ s_i \in S$ and $a_i \overline{s}_i\prec w$, then this
composition is called trivial modulo $(S,w)$ and denote it by
$$
(f,g)_w\equiv 0 \ mod(S,w).
$$

\begin{definition}\label{d3}(\cite{ch})
Let $S\subset mod_{k\langle X\rangle} \langle Y\rangle$ be a
non-empty set  with each $s\in S$ monic. Let the order ``$\prec$" be
as before. Then we call $S$ a Gr\"obner-Shirshov basis in the module
$mod_{k\langle X\rangle} \langle Y\rangle$ if all the compositions
in $S$ are trivial modulo $S$.
\end{definition}

The proof of the following lemma is essentially from \cite{ch}. For
convenience, we give the detail.

\begin{lemma}\label{l3.5}
(\cite{ch}, Composition-Diamond lemma for ``double-free"
modules)\label{l2} Let $S\subset mod_{k\langle X\rangle} \langle
Y\rangle$ be a non-empty set  with each $s\in S$ monic and
$``\prec"$ the order on $X^*Y$ as before. Then the following
statements are equivalent:
\begin{enumerate}
\item[(1)]\ $S$ is a Gr\"obner-Shirshov basis in $mod_{k\langle X\rangle} \langle Y\rangle$.
 \item[(2)]\ If \  $0\neq
f\in k\langle X\rangle S$, then $\overline{f}=a\overline{s}$ for
some $a\in X^*, \ s\in S$.
\item[(2')]\ If \ $0\neq f\in k\langle
X\rangle S$, then $f=\sum \alpha_i a_i s_i$ with $a_1
\overline{s}_1\succ a_2 \overline{s}_2\succ\cdots$, where each
$\alpha_i\in k, \ a_i\in X^*, \ s_i\in S$.
\item[(3)]\
$Red(S)=\{w \in X^*Y|w\neq a\overline{s}, \ a\in X^*, \ s\in S\}$ is
a $k$-linear basis for the factor $mod_{k\langle X\rangle} \langle
Y|S\rangle=mod_{k\langle X\rangle} \langle Y\rangle /k\langle
X\rangle S$.
\end{enumerate}
\end{lemma}
\textbf{Proof: }$(1)\Rightarrow (2)$. Suppose that $0\neq f\in
k\langle X\rangle S$. Then $f=\sum \alpha_i a_i s_i$ for some
$\alpha_i\in k, \ a_i\in X^*, \ s_i\in S$. Let $w_i=a_i
\overline{s}_i$ and $w_1=w_2=\cdots=w_l\succ w_{l+1}\succeq \cdots$.
We will prove that $\overline{f}=a\overline{s}$ for some $a\in X^*,
\ s\in S$, by using induction on $l$ and $w_1$. If $l=1$, then the
result is clear. If $l>1$, then $a_1 \overline{s}_1=a_2
\overline{s}_2$. Thus, we may assume that $a_1=a_2 a, \
\overline{s}_2=a\overline{s}_1$ for some $a\in X^*$. Now, by (1),
$$
a_1 s_1-a_2 s_2=a_2 as_1-a_2 s_2=a_2(as_1-s_2)=a_2\sum \beta_j b_j
u_j=\sum \beta_j a_2b_j u_j,
$$
where $\beta_j\in k, \ b_j\in X^*, \ u_j\in S$ and $b_j
\overline{u}_j\prec\overline{s}_2$. Therefore, $a_2 b_j
\overline{u}_j\prec w_1$. Now, by induction on $l$ and $w_1$, we
have the result.

It is clear that $(2)$ is equivalent to $(2')$.

$(2)\Rightarrow (3)$. For any $0\neq f\in mod_{k\langle X\rangle}
\langle Y\rangle$, if $\bar{f}=u_1\in Red(S)$, then $f=\beta_1
u_1+\cdots$. If $\bar{f}\not\in Red(S)$, then $f=\alpha_1 a_1
s_1+\cdots$. It follows that we can express $f$ as
$$
f=\sum \alpha_i a_i s_i+\sum \beta_j u_j,
$$
where $\alpha_i, \ \beta_j \in k, \ a_i\in X^*, \ s_i\in S \
\mbox{and} \ u_j\in Red(S)$. Then $Red(S)$ generates the factor
module. Moreover, assume that $0\neq \sum \alpha_i a_i s_i=\sum
\beta_j u_j$, where $ \ a_i\in X^*, \ s_i\in S, \ u_j\in Red(S), \
a_1 \overline{s}_1\succ a_2 \overline{s}_2\succ\cdots$ and $u_1\succ
u_2\succ\cdots$. Then $u_1=a_1 \overline{s}_1$, a contradiction.
This shows that $Red(S)$ is a $k$-linear basis of the factor module.

$(3)\Rightarrow (1)$. For any $f, \ g\in S$, suppose that
$w=\overline{f}=a\overline{g}$. Since $(f,g)_w\in k\langle X\rangle
S$,  we get, by $(3)$, that $(f,g)_w=f-ag=\sum \alpha_i a_i s_i$,
where $s_i\in S, \  a_i\in X^*$ and $a_i \overline{s_i}\preceq
\overline{(f,g)_w}\prec w$.  Then, $S$ is a Gr\"obner-Shirshov basis
in the module $mod_{k\langle X\rangle} \langle Y\rangle$. \ \
$\square$
 \\

\noindent{\bf Remark}: Let us view $k\langle X\rangle$ as free left
$k\langle X\rangle$-module with one generator $I$. Then
$mod_{k\langle X\rangle} \langle I\rangle=k\langle X\rangle
I=_{k\langle X\rangle}\!\!k\langle X\rangle$ is a cyclic $k\langle
X\rangle$-module. If $S\subset k\langle X\rangle$, then $k\langle
X\rangle S$ is a left ideal of $k\langle X\rangle$ (also a left
$k\langle X\rangle$-submodule of $k\langle X\rangle I$). This
implies that Lemma \ref{l1} is a special case of Lemma \ref{l2}.

\ \

Let $S\subset k\langle X\rangle$ and $A=k\langle X|S\rangle$ an
associative algebra.  Then, for any left $A$-module $_AM$, we can
view $_AM$ as a $k\langle X\rangle$-module in a natural way: for any
$f\in k\langle X\rangle, \ m\in M$,
$$
fm=(f+Id(S))m.
$$
We note that $_AM$ is an epimorphic image of some free $A$-module.
Then, we can assume that $_A M=mod_A \langle Y|T\rangle=mod_{A}
\langle Y\rangle/AT$, where $T\subset mod_{A} \langle Y\rangle$ and
$mod_{A} \langle Y\rangle$ a free left $A$-module with the basis
$Y$. Let  $T_1=\{\sum f_iy_i\in mod_{k\langle X\rangle} \langle
Y\rangle|\sum (f_i+Id(S))y_i\in T\}$ and $R=SX^*Y\cup T_1$. Then, by
the following Lemma \ref{l3.6}, we have, as $k\langle
X\rangle$-modules, $_A M\cong mod _{k\langle X\rangle}\!\langle
Y|R\rangle$.

\begin{lemma}(\cite{ch})\label{l3.6} Let the notations be as above.
Then, as $k\langle X\rangle$-modules,
$$
\sigma: \ _A M\rightarrow mod _{k\langle X\rangle}\!\langle
Y|R\rangle, \ \ \sum (f_i+Id(S)) (y_i+AT)\mapsto\sum f_i
y_i+k\langle X\rangle R
$$
is an isomorphism, where each $f_i\in k\langle X\rangle$.
\end{lemma}

\textbf{Proof: } For any $\sum (f_i+Id(S)) (y_i+AT),\sum (g_i+Id(S))
(y_i+AT)\in_A\!M$, we have
\begin{eqnarray*}
&&\sum (f_i+Id(S)) (y_i+AT)=\sum (g_i+Id(S)) (y_i+AT)\ \ \ \mbox{ in } \ \ \ _A\!M \\
&\Leftrightarrow& \sum (f_i-g_i)y_i\in AT\ \ \ \mbox{ in } \ \ \ _A\!M \\
&\Leftrightarrow&  \sum (f_i-g_i)y_i\in k\langle X\rangle R\\
&\Leftrightarrow&  \sum f_i y_i+k\langle X\rangle R=\sum g_i
y_i+k\langle X\rangle R.
\end{eqnarray*}
Hence, $\sigma$ is injective. It is easy to see that $\sigma$ is a
surjective mapping and a $k\langle X\rangle$-module homomorphism. \
\ $\square$

By using Lemma \ref{l3.5} and Lemma \ref{l3.6}, we know that if we
want to find a $k$-linear basis for the module $_A M=mod_A \langle
Y|T\rangle$, where $A=k\langle X|S\rangle$, we only need to find a
Gr\"obner-Shirshov basis for the module $mod_{k\langle X\rangle}
\langle Y|SX^*Y\cup T_1\rangle$, where $T_1=\{\sum f_iy_i\in
mod_{k\langle X\rangle} \langle Y\rangle|\sum (f_i+Id(S))y_i\in
T\}$.

\ \

The following theorem gives some relationships between  the
Gr\"{o}bner-Shirshov bases (pairs) in free associative algebras and
``double-free" modules.

\begin{theorem}\label{l3}
Let $X, \ Y$ be well ordered sets, $<$ a monomial order on $X^*$,
$\prec$ the order on $X^*Y$ as in $(*)$. Let $S,\ T\subset k\langle
X\rangle$ be monic sets. Then the following statements hold:
\begin{enumerate}
\item[(1)]\ $S\subset k\langle X\rangle$
is a Gr\"obner-Shirshov basis in $k\langle X\rangle$  with respect
to the order $<$ if and only if $SX^*Y\subset mod_{k\langle
X\rangle} \langle Y\rangle$ is a Gr\"obner-Shirshov basis in
$mod_{k\langle X\rangle} \langle Y\rangle$ with respect to the order
$\prec$.
 \item[(2)]\ Let us view $k\langle X\rangle$ as
free $k\langle X\rangle$-module with one generator $I$. Then,
$(S,T)$ is a Gr\"obner-Shirshov pair for the $A$-module
$M=A/A(T+Id(S))$, where $A=k\langle X|S\rangle$, if and only if $S$
is a Gr\"obner-Shirshov basis in the algebra $k\langle X\rangle$
with respect to the order $<$ and $(SX^*\cup T)I$ a
Gr\"obner-Shirshov basis in the free module $mod_{k\langle
X\rangle}\!\langle I\rangle$ with respect to the order $\prec$.
\end{enumerate}
\end{theorem}
\textbf{Proof: } (1) Suppose that $S$ is a Gr\"obner-Shirshov basis
in $k\langle X\rangle$. We shall prove that all the compositions in
$SX^*Y$ are trivial modulo $SX^*Y$. For any $f, \ g\in SX^*Y$, let
$f=s_1 a_1 y, \ g=s_2 a_2 y , \ s_1,s_2\in S, \ a_1,a_2\in X^*, \
y\in Y$ and $w=\overline{f}=a\overline{g}$. Then $\overline{s}_1
a_1=a\overline{s}_2 a_2$. Since $S$ is a Gr\"obner-Shirshov basis in
$k\langle X\rangle$, we have
$$
(f,g)_w=f-ag=s_1 a_1 y-as_2 a_2 y=(s_1 a_1-as_2 a_2)y=\sum (\alpha_i
u_i r_i v_i)y,
$$
where $u_i, \ v_i\in X^*, \ r_i\in S$ and $u_i \overline{r}_i v_i
y\prec w$. So each composition is trivial modulo $SX^*Y$ and hence,
$SX^*Y$ is a Gr\"obner-Shirshov basis in  $mod_{k\langle X\rangle}
\langle Y\rangle$.

Conversely, assume that $SX^*Y$ is a Gr\"obner-Shirshov basis in the
module $mod_{k\langle X\rangle} \langle Y\rangle$. For any $f, \
g\in S$ and $w=\overline{f}a=b\overline{g}$, we have
$w_1=\overline{fay}=b\overline{gy}$ and
$$
(fay,bgy)_{w_1}=(fa-bg)y=\sum \alpha_i (a_i r_i) y,
$$
where $\alpha_i\in k, \ r_i=s_i b_i, \ a_i,b_i\in X^*, \ s_i \in S$
and $a_i \overline{r}_i y\prec w_1$. Then,
$$
(f,g)_w=fa-bg=\sum \alpha_i a_is_i b_i
$$
with $a_i \overline{s}_i b_i<w$. This shows that, each composition
of intersection in $S$ is trivial modulo $S$. Similarly, each
composition of inclusion in $S$ is trivial modulo $S$. Therefore,
$S$ is a Gr\"obner-Shirshov basis in $k\langle X\rangle$.

(2)  The results follow from  the Definitions \ref{d4} and \ref{d3}
directly.
 \ \ $\square$

 \ \

\noindent{\bf Remark}: By Theorem \ref{l3}, it is clear that the
Theorem \ref{t1} follows from Lemma \ref{l3.5}.

\section{Highest weight
modules over $sl_2$}

In this section, as an application of Lemma \ref{l2}, we re-prove
that the highest weight module over $sl_2$ is irreducible (see
\cite{book l}) and we show that any finite dimensional irreducible
$sl_2$-module has the presentation $(**)$ as below.

Let $X=\{x, \ y, \ h\}$ and $sl_2=Lie(X| S)$ a Lie algebra over a
field $k$ with $chk=0$, where
\[
x=\left(
\begin{array}{ll}
0 & 1 \\
0 & 0
\end{array}
\right), \ \ \
y=\left(
\begin{array}{ll}
0 & 0 \\
1 & 0
\end{array}
\right), \ \ \
h=\left(
\begin{array}{ll}
1 &  0 \\
0 & -1
\end{array}
\right) \ \mbox{ and } S=\{[hx]-2x, \ [hy]+2y, \ [xy]-h\}.
\]
Then the universal enveloping algebra of $sl_2$ is $\mathcal
{U}(sl_2)=k\langle X|\ S^{(-)}\rangle$. Define the deg-lex order on
$X^*$ with $x>h>y$. Then $S$ is a Gr\"obner-Shirshov basis in free
Lie algebra $Lie(X)$ since $S^{(-)}$ is a Gr\"obner-Shirshov basis
in $k\langle X\rangle $ (see Lemma \ref{l0}). Let
$$_{sl_2} V(\lambda)=mod_{sl_2} \langle v_0| xv_0=0, \
hv_0=\lambda v_0, \ y^{m+1}v_0=0\rangle
$$
be a highest weight module generated by $v_0$ with the highest
weight $\lambda$. We can rewrite it as
\begin{eqnarray*}
_{sl_2} V(\lambda)&=&mod_{\mathcal {U}(sl_2)} \langle v_0| \ xv_0=0,
\ hv_0=\lambda v_0, \ y^{m+1}v_0=0 \rangle\\
&=&mod_{k\langle X\rangle} \langle v_0| \ xv_0=0, \ hv_0=\lambda
v_0, \ y^{m+1}v_0=0, \ S^{(-)}X^*v_0=0\rangle
\end{eqnarray*}
Let $S_1=\{xv_0, \ hv_0-\lambda v_0, \ y^{m+1}v_0\}\cup
S^{(-)}X^*v_0$. It is easy to see that all compositions in $S_1$ are
trivial modulo $S_1$. Thus, $S_1$ is a Gr\"obner-Shirshov basis for
this module with respect to the order as $(*)$ in \S 3, and by Lemma
\ref{l2}, $Red(S_1)=\{y^i v_0|0\leq i\leq m\}$ is a $k$-linear basis
for module $_{sl_2} V(\lambda)$, and so $dim(V(\lambda))=m+1$.

Let $y^{(i)}=\frac{1}{i!}y^i, \ v_i=\frac{1}{i!}y^iv_0$ and
$v_{-1}=0$. Then $v_i \ (0\leq i\leq m)$ is also a linear basis of
$V(\lambda)$. Now, using $ELW$ of the relations in $S_1$ on the left
parts, we have the following equalities (see also \cite{book l},
p.32):
\begin{lemma}\label{l4}
\begin{eqnarray*}
hv_i&=&(\lambda-2i)v_i\\
yv_i&=&(i+1)v_{i+1}\\
xv_i&=&(\lambda-i+1)v_{i-1} \ (0\leq i) \ \ \ \ \ \square
\end{eqnarray*}
\end{lemma}

Since $v_{m+1}=0$ and $chk=0$, $0=xv_{m+1}=(\lambda-m)v_m$ and
therefore, $\lambda=m$.

\begin{lemma}
$V(\lambda)$ is irreducible.
\end{lemma}
\textbf{Proof: }Let $0\neq V_1\leq V(\lambda)$ be a submodule. Since
$V_1\neq0$, there exist $0\neq a_i v_i+a_{i+1} v_{i+1}+\cdots+ a_m
v_m$, where $i$ is the least number such that $a_i\neq 0$. Applying
$y$ to it $m-i$ times, we get $a_i(i+1)(i+2)\cdots m v_m\in V_1$.
So, $v_m\in V_1$. Applying $x$ to $v_m$, we get $v_i\in V_1 \ (0\leq
i<m)$ and $V_1= V(\lambda)$.
 \ \ $\square$

For any finite dimensional irreducible $sl_2$-module $V$, choosing a
maximal vector $v_0\in V$ and $v_i=\frac{1}{i!}y^i v_0$, we have the
formulas as in Lemma \ref{l4}. We can suppose that $dimV=m$. Thus,
$v_m\neq 0, \ v_{m+1}=0$ and hence, $V$ can be represented as
$$
_{sl_2} V=mod_{sl_2} \langle v_0| xv_0=0, \ hv_0=\lambda v_0, \
y^{m+1}v_0=0\rangle \ \ \ \ \ \ \  \ \ \ \ \ \ \ \ \ \ \ \ \ (**)
$$
which means that any finite dimensional irreducible $sl_2$-module
has the above form.

\section{Verma module over Kac-Moody algebras}

In this section, we give the definitions of Kac-Moody algebra
$\mathcal {G}(A)$ and Verma module over $\mathcal {G}(A)$. Then, by
using Lemma \ref{l3.5}, we find a Gr\"obner-Shirshov basis for this
Verma module.

Let $A=(a_{ij})$ be an (integral) symmetrizable n-by-n Cartan matrix
over $\mathbb{C}$, where $\mathbb{C}$ is the complex field. It means
that $a_{ii}=2, \ a_{ij}\leq 0 \ (i\neq j)$, and there exists a
diagonal matrix $D$ with diagonal entries $d_i$ nonzero integers
such that product $DA$ is symmetric. Let $\mathcal {G}(A)=Lie(Z|S)$
be a Lie algebra, where $Z=\{x_i, \ y_i, \ h|1\leq i \leq n, \ h\in
H \}$, $S$ consists of the following relations (see \cite{book 2},
p.159):
\begin{enumerate}
\item[(5.1)]\ $[x_i, y_j]=\delta_{ij}\alpha_i^{\vee } \ \ \
(i,j=1,\cdots,n)$,
\item[(5.2)]\ $[h,h']=0  \ \ \  \ \ (h, h'\in H)$,
\item[(5.3)]\ $[h,x_i]=\langle \alpha_i, h\rangle x_i,\ [h,y_i]=-\langle \alpha_i, h\rangle y_i ,\ (i=1,\cdots,n; h\in
H)$,
\item[(5.4)]\ $(adx_i)^{1-a_{ij}}x_j=0,\ (ady_i)^{1-a_{ij}}y_j=0 \ \ \ \ \ (i\neq
j)$,
\end{enumerate}
where $ad$ is the derivation, $H$ a complex vector space,
$\Pi=\{\alpha_1, \cdots,\alpha_n\}\subset H^{\star}$ (the dual space
of H) and $\Pi^{\vee }=\{\alpha_1^{\vee }, \cdots,\alpha_n^{\vee
}\}\subset H$ indexed subsets in $H^{\star}$ and $H$, respectively,
satisfying the following conditions (see \cite{book 2}, p.1):

\begin{enumerate}
\item[(a)]\ both sets $\Pi$ and $\Pi^{\vee }$ are linearly
independent, \item[(b)]\ $\langle\alpha_i^{\vee }, \alpha_j
\rangle=a_{ij} \ \ \ \ \ (i, \ j=1,\cdots,n)$, \item[(c)]\
$n-l=dimH-n \ \ \ \ \ rank(A)=l.$
\end{enumerate}

Then we call this Lie algebra $\mathcal {G}(A)$ Kac-Moody algebra.
Let $\mathfrak{N}_+ \ \  (\mathfrak{N}_-)$ be the subalgebra of
$\mathcal {G}(A)$ generated by $x_i \ \ (y_i), \ (0\leq i\leq n)$.
Then $\mathcal {G}(A)=\mathfrak{N}_-\oplus H \oplus \mathfrak{N}_+$
and $\mathcal {U}(\mathcal {G}(A))=\mathcal
{U}(\mathfrak{N}_+)\otimes k[H]\otimes \mathcal {U}(\mathfrak{N}_-)$
is the universal enveloping algebra of $\mathcal {G}(A)$, where
$\mathcal {U}(\mathfrak{N}_+) \ (\mathcal {U}(\mathfrak{N}_-))$ is
the universal enveloping algebra of $\mathfrak{N}_+ \
(\mathfrak{N}_-)$.

Let $\{h_j | 1\leq j\leq 2n-l\}$ be a basis of $H$. We order the set
$X=\{x_i, \ h_j, \ y_m|1\leq i,m \leq n, \ 1\leq j\leq 2n-l\}$ by
$x_i>x_j, \ h_i>h_j, \
 y_i>y_j$, if $i>j$, and $x_i>h_j>y_m$ for all $i, \ j, \ m$.
Then we define the deg-lex order on $ X^*$.

By \cite{bm96}, we can get a Gr\"obner-Shirshov basis $T$ for
$\mathcal {U}(\mathcal {G}(A))$, where $T$ consists of the following
relations:

\begin{enumerate}
\item[(5.5)]\ $h_ih_j-h_jh_i, \ x_jh_i-h_ix_j+d_ia_{ij}x_i, \
h_iy_j-y_jh_i+d_ia_{ij}y_j$,
\item[(5.6)]\ $x_iy_j-y_jx_i-\delta_{ij}h_i$,
\item[(5.7)]\ $\{\sum\limits_{\nu=0}^{1-a_{ij}} (-1)^\nu \left[
\begin{array}{ll}
1-a_{ij}\\
\nu
\end{array}
\right]x_i^{1-a_{ij}-\nu}x_jx_i^{\nu}\}^c  \ \ \ \ \ \ \ \ (i\neq
j)$,
\item[(5.8)]\ $\{\sum\limits_{\nu=0}^{1-a_{ij}} (-1)^\nu \left[
\begin{array}{ll}
1-a_{ij}\\
\nu
\end{array}
\right]y_i^{1-a_{ij}-\nu}y_jy_i^{\nu}\}^c \ \ \ \ \ \ \ \ (i\neq
j).$
\end{enumerate}
where $S^c$ means a Gr\"obner-Shirshov basis that contains $S$.

\begin{definition}\em(\cite{book 2}) A $\mathcal {G}(A)$-module $V$
is called a highest weight module with highest weight $\Lambda\in
H^{\star}$ if there exist a non-zero vector $v\in V$, such that
$$
\mathfrak{N}_+(v)=0, \ h(v)=\Lambda(h)v, \ h\in H
$$
and $\mathcal {U}(\mathcal {G}(A))(v)=V$.
\end{definition}

A Verma module $M(\Lambda)$ with highest weight $\Lambda$ has the
following presentation:
$$
_{\mathcal {G}(A)}M(\Lambda)=mod_{\mathcal {G}(A)}\langle
v|\mathfrak{N}_+(v)=0, \ h(v)=\Lambda(h)v, \  \ h\in H\rangle.
$$

\begin{corollary}
$R=\{TX^*(v), \ \mathfrak{N}_+(v), \ h(v)-\Lambda(h)v\}$ is a
Gr\"obner-Shirshov basis for the Verma module $_{\mathcal
{G}(A)}M(\Lambda)$.
\end{corollary}
\textbf{Proof: } Since
\begin{eqnarray*}
_{\mathcal {G}(A)}M(\Lambda)&=&mod_{\mathcal {U}(\mathcal
{G}(A))}\langle v|\mathfrak{N}_+(v)=0, \ h(v)=\Lambda(h)v, \  \ h\in H\rangle\\
&=&mod _{k\langle Z\rangle}\langle v|TX^*(v)=0, \
\mathfrak{N}_+(v)=0, \ h(v)=\Lambda(h)v, \ h\in H\rangle,
\end{eqnarray*}
it is easily to check that all compositions in $R$ are trivial. So,
$R$ is Gr\"obner-Shirshov basis for the Verma module. \ \ $\square$

\noindent{\bf Remark}: In the book \cite{book l}, the author only
consider the semisimple Lie algebras and call this highest weight
module to be standard cyclic module.

\section{Verma module over the coefficient algebra of a free Lie conformal algebra}
In this section, by using Lemma \ref{l3.5}, we find a basis of Verma
module over Lie algebras of coefficients of free conformal algebras.

Let $\mathcal {B}$ be a set of symbols. Let the locality function
$N: \mathcal {B}\times \mathcal {B}\rightarrow \mathbb{Z}_+$ be a
constant, i.e. $N(a,b)\equiv N$ for any $a, \ b\in \mathcal {B}$.
Let $X=\{b(n)| \ b\in \mathcal {B}, \ n\in \mathbb{Z}\}$ and
$L=Lie(X|S)$ be a Lie algebra generated by $X$ with the relation
$S$, where
\[
S=\{\sum\limits_{s} (-1)^s \left(
\begin{array}{ll}
n \\
s
\end{array}
\right)
 [b(n-s)a(m+s)]=0| \ a, \ b\in \mathcal {B}, \ m, n\in
\mathbb{Z}\}.
\]
For any $b\in \mathcal {B}$, let $\widetilde{b}=\sum\limits_{n} b(n)
z^{-n-1}\in L[[z,z^{-1}]]$. It is well-known that they generate a
free Lie conformal algebra $C$ with data $(\mathcal {B},N)$ (see
\cite{ro99}). Moreover, the coefficient algebra of $C$ is just $L$.

Let $\mathcal {B}$ be a linearly ordered set. Define an order on $X$
in the following way:
$$
a(m)<b(n)\Leftrightarrow m<n \ \mbox{ or } (m=n \ \mbox{and} \ a<b).
$$
We use the deg-lex order on $X^*$. Then, it is clear that the
leading term of each polynomial in $S$ is $b(n)a(m)$ such that
$$
n-m>N \ \mbox{or} \ (n-m=N \ \mbox{and} \ (b>a \ \mbox{or} \ (b=a \
\mbox{and} \ N \ \mbox{is odd})).
$$

The following lemma is  from \cite{ro99}.

\begin{lemma}
$S$ is a Gr\"obner-Shirshov basis in $Lie(X)$.
\end{lemma}

\begin{corollary}
Let $\mathcal {U}=\mathcal {U}(L)$ be the universal enveloping
algebra of $L$. Then a $k$-basis of $\mathcal {U}$  consists of
monomials
$$a_1(n_1)a_2(n_2)\cdots a_k(n_k), \ a_i\in \mathcal {B}, \ n_i\in \mathbb{Z}$$
such that for any $1\leq i<k$,
\[
n_i-n_{i+1}\leq \left \{
\begin{array}{ll}
N-1 & \mbox{if $a_i>a_{i+1}$ or ($a_i=a_{i+1}$ \ and \ $N$ is odd)} \\
 N  & \mbox{otherwise}
 \end{array}
 \right. \ \ \ \ \ \ \ \  \ \ \ \  \ (***)
 \]
 \end{corollary}

\textbf{Proof: } Viewing $\mathcal {U}$ as a $k\langle
X\rangle$-module, we have
$$
_{\mathcal {U}} \mathcal {U}=mod_{\mathcal {U}}\langle X| \
S^{(-)}\rangle=mod_{k\langle X\rangle} \langle I|
 \ S^{(-)}X^*I\rangle.
$$
 Since $S$ is a Gr\"obner-Shirshov basis
in $Lie(X)$, $S^{(-)}$ is a Gr\"obner-Shirshov basis in $k\langle
X\rangle$ by Lemma \ref{l0}. Therefore, by Theorem \ref{l3},
$S^{(-)}X^*I$ is a Gr\"obner-Shirshov basis in the free module
$mod_{k\langle X\rangle}\!\langle I\rangle$. Now, the result follows
from Lemma \ref{l3.5}. \ \ $\square$

\begin{definition}\em (\cite{book 2} \cite{ka97})
\begin{enumerate}
\item[(a)]\ An $L$-module $M$ is called restricted, if for any $a\in C, \
v\in M$ there is some integer $T$ such that for any $n\geq T$ one
has $a(n)v=0$.
\item[(b)]\ An $L$-module $M$ is called highest weight module
if it is generated over $L$ by a single element $m\in M$ such that
$L_+ m=0$, where $L_+$ is the subspace of $L$ generated by
$\{a(n)|a\in C, n\geq0\}.$ In this case $m$ is called the highest
weight vector.
\end{enumerate}
\end{definition}

Now we build a universal highest weight module $V$ over $L$, which
is often referred to as Verma module. Let $kI_v$ be a
$1$-dimensional trivial $L_+$-module generated by $I_v$, i.e.,
$a(n)I_v=0$ for all $a\in \mathcal {B}, \ n\geq0$. Clearly,
$$
V=Ind_{L_+} ^L kI_v=\mathcal {U}(L)\otimes_{\mathcal {U}(L_+)}
kI_v\cong \mathcal {U}(L)/\mathcal {U}(L)L_+.
$$
Then $V$ has a structure highest weight module over $L$ with the
action given by the multiplication on $\mathcal {U}(L)/\mathcal
{U}(L)L_+$ and the highest weight vector $I\in \mathcal {U}(L)$.
Also $V=\mathcal {U}(L)/\mathcal {U}(L)L_+$ is the universal
enveloping vertex algebra of $C$ and the embedding $\varphi:
C\rightarrow V$ is given by $a\mapsto a(-1)I$ (see also
\cite{ro99}).

\begin{theorem}
Let the notions be as above. Then a $k$-basis of $V$ consists of
elements
$$
a_1(n_1)a_2(n_2)\cdots a_k(n_k), \ a_i\in \mathcal {B}, \ n_i\in
\mathbb{Z}
$$
such that the condition $(***)$ holds and $n_k<0$.
\end{theorem}

\textbf{Proof: } Clearly, as $k\langle X\rangle$-modules,
$$
_{\mathcal {U}} V=_{\mathcal {U}} (\mathcal {U}(L)/\mathcal
{U}(L)L_+)=mod_{k\langle X\rangle} \langle I| \ S^{(-)}X^*I, \
a(n)I, \ n\geq0\rangle=_{k\langle X\rangle} \langle I| \ S'\rangle,
$$
where $S'=\{S^{(-)}X^*I, \ a(n)I, \ n\geq0\}$.  In order to prove
that $S'$ is a Gr\"obner-Shirshov basis, we only need to check
$w=b(n)a(m)I$, where $m\geq 0$. Let
\[
f=\sum\limits_{s} (-1)^s \left(
\begin{array}{ll}
n \\
s
\end{array}
\right)
 (b(n-s)a(m+s)-a(m+s)b(n-s))I \ \ and \ \ \ \
 g=a(m)I.
\]
Then $(f,g)_w=f-b(n)a(m)I\equiv 0 \ mod(S',w)$ since $n-m\geq N$,
$m+s\geq0, \ n-s\geq0, \ 0\leq s \leq N$. It follows that $S'$ is a
Gr\"obner-Shirshov basis.  Now, the result follows from Lemma
\ref{l3.5}. \ \ $\square$
\\

\section{Universal enveloping module for a
Sabinin algebra}

In this section, we deal with the universal enveloping module for a
Sabinin algebra.
\begin{definition}\em(\cite{pe})
A vector space $V$ is called a Sabinin algebra if it is endowed with
multilinear operation $\langle;\rangle$: for any
$x_1,x_2,\cdots,x_m,y,z\in V$ and any $m\geq0$,
$$
 \langle x_1,x_2,\cdots,x_m;y,z\rangle
$$
satisfying the identities
\begin{eqnarray*}
&& \langle x_1,x_2,\cdots,x_m;y,
z\rangle =-\langle x_1,x_2,\cdots,x_m;z,y\rangle,\\
&& \langle x_1,x_2,\cdots, x_r,a,b,x_{r+1},\cdots,x_m; y,z\rangle
-\langle x_1,x_2,\cdots, x_r,b,a,x_{r+1},\cdots,x_m; y,z\rangle\\
&& \ \  \ \ \ +\sum_ {k=0}^ {r}\sum_ {\alpha}\langle
x_{\alpha_1},\cdots,x_{\alpha_k},\langle
x_{\alpha_{k+1}},\cdots ,x_{\alpha_r};a,b\rangle,\cdots,x_m;y,z\rangle =0,\\
&& \sigma_{x,y,z}(\langle x_1,x_2,\cdots, x_r,x; y,z\rangle+\sum_
{k=0}^ {r}\sum_ {\alpha}\langle
x_{\alpha_1},\cdots,x_{\alpha_k};\langle x_{\alpha_{k+1}},\cdots
,x_{\alpha_r};y,z\rangle,x\rangle) =0,
\end{eqnarray*}
where $\alpha$ runs the set of all bijections of the type $\alpha:\
\{1,2,\cdots,r\}\rightarrow \{1,2,\cdots,r\},\ i\mapsto \alpha_i,\
\alpha_1<\alpha_2<\cdots<\alpha_k,\ \alpha_{k+1}<\cdots<\alpha_r,\
r\geq0,\ \sigma_{x,y,z}$ denotes the cyclic sum by $x,y,z$.
\end{definition}

Let $X=\{a_i|i\in \Lambda\}$ be a totally ordered basis of $V$. We
define the deg-lex order on $X^*$. Let $\Delta: \ V\rightarrow
V\otimes V$ be a linear map which satisfies: $\Delta(a_i)=1\otimes
a_i+a_i\otimes 1, \ (Id\otimes\Delta)\Delta=(\Delta\otimes Id)\Delta
\mbox{ (coassociative)}\ \mbox{\ and \ if }\  \tau\Delta=\Delta \
\mbox{\ then }\  \tau(x\otimes y)=y\otimes x \ \mbox{
(cocommutative)}$. It is customary to write $\Delta(x)=\sum
x_{(1)}\otimes x_{(2)}$.

Let $T(V)$ be the tensor algebra over $V$ endowed with its natural
structure of cocommutative Hopf algebra, that is, $V\subseteq
Prim(T(V))$ (the primitive element of $T(V)$). Let $\langle;\rangle:
T(V)\bigotimes V \bigotimes V\rightarrow V$ be a map. Then we may
write the definition of a Sabinin algebra very shortly as
\begin{eqnarray*}
&&\ \ \ \ \ \ \ \ \ \ \ \ \ \ \ \ \langle x;a,
b\rangle =-\langle x;b,a\rangle,\\
&&\\
 &&\ \ \ \ \ \langle x[a,b]y;
c,e\rangle
+\sum\langle x_{(1)}\langle x_{(2)};a,b\rangle y; c,e\rangle=0,\\
&&\\
&&\ \ \ \ \ \sigma_{a,b,c}(\langle xc; a,b\rangle+\sum\langle
x_{(1)};\langle x_{(2)};a,b\rangle,c\rangle) =0,
\end{eqnarray*}
where $[a,b]=ab-ba$.

\begin{definition}\em(\cite{pe}) Let $(V,\langle;\rangle)$ be a Sabinin
algebra. Then
$$
\widetilde{S}(V)=T(V)/span\langle xaby-xbay+\sum x_{(1)}\langle
x_{(2)};a,b\rangle y|x,y\in T(V),a,b\in V\rangle
$$
is called the universal enveloping module for $V$.
\end{definition}

 Since $T(V)\simeq k\langle X\rangle$ as $k$-algebras, we
can view $\widetilde{S}(V)$ as a right $k\langle X\rangle$-module:
$$
\widetilde{S}(V)=mod\langle X|I\rangle_{k\langle X\rangle},
$$
where $ I=\{xab-xba+\sum x_{(1)}\langle x_{(2)};a,b\rangle|x\in
X^*,a>b, \ a,b\in X\}. $

Then we have the following theorem.
\begin{theorem}\label{t7.1}
Let $I$ be as above. Then $I$ is the Gr\"obner-Shirshov basis in
$mod\langle X\rangle_{k\langle X\rangle}$.
\end{theorem}
\textbf{Proof: } There are two kinds of compositions: $w_1=xabc\
(a>b>c) \ \mbox{ and } \ w_2=ucdvab\ (c>d,a>b)$. Denote by
\begin{eqnarray*}
f_1&=&xabc-xacb+\sum(xa)_{(1)}\langle (xa)_{(2)};b,c\rangle,\\
f_2&=&xab-xba+\sum x_{(1)}\langle x_{(2)};a,b\rangle,\\
f_3&=&ucdvab-ucdvba+\sum(ucdv)_{(1)}\langle
(ucdv)_{(2)};a,b\rangle,\\
f_4&=&ucd-udc+\sum u_{(1)}\langle u_{(2)};c,d\rangle.
\end{eqnarray*}
Then
\begin{eqnarray*}
(f_1,f_2)_{w_1}&=&xabc-xacb+\sum x_{(1)}a\langle
x_{(2)};b,c\rangle+\sum
x_{(1)}\langle x_{(2)}a;b,c\rangle\\
&&-xabc+xbac-\sum x_{(1)}\langle
x_{(2)};a,b\rangle c\\
&\equiv&-xcab+\sum x_{(1)}\langle x_{(2)};a,c\rangle b+\sum
x_{(1)}a\langle
x_{(2)};b,c\rangle+\sum x_{(1)}\langle x_{(2)}a;b,c\rangle\\
&&+xbca-\sum x_{(1)}b\langle x_{(2)};a,c\rangle -\sum x_{(1)}\langle
x_{(2)}b;a,c\rangle-\sum x_{(1)}\langle x_{(2)};a,b\rangle c\\
 &\equiv&\sum x_{(1)}c\langle x_{(2)};a,b\rangle+\sum
x_{(1)}\langle
x_{(2)}c;a,b\rangle+\sum x_{(1)}\langle x_{(2)};a,c\rangle b\\
&&+\sum x_{(1)}a\langle x_{(2)};b,c\rangle+\sum x_{(1)}\langle
x_{(2)}a;b,c\rangle -\sum x_{(1)}\langle x_{(2)};b,c\rangle a\\
&&-\sum x_{(1)}b\langle x_{(2)};a,c\rangle-\sum x_{(1)}\langle
x_{(2)}b;a,c\rangle-\sum x_{(1)}\langle x_{(2)};a,b\rangle c\\
&\equiv&\sum x_{(1)}\langle x_{(2)}a;b,c\rangle+\sum x_{(1)}\langle
x_{(2)};\langle x_{(3)};b,c\rangle,
a\rangle+\sum x_{(1)}\langle x_{(2)}b;c,a\rangle\\
&&+\sum x_{(1)}\langle x_{(2)}c;a,b\rangle +\sum x_{(1)}\langle
x_{(2)};\langle x_{(3)};c,a\rangle, b\rangle \\
&&+\sum x_{(1)}\langle x_{(2)};\langle x_{(3)};a,b\rangle,
c\rangle\\
&\equiv&0
\end{eqnarray*}
since $\sigma_{a,b,c}(\langle xc; a,b\rangle+\sum\langle
x_{(1)};\langle x_{(2)};a,b\rangle,c\rangle) =0$.

\begin{eqnarray*}
(f_3,f_4)_{w_2}&=& ucdvab-ucdvba+\sum u_{(1)}v_{(1)}\langle
u_{(2)}cdv_{(2)};a,b\rangle+\sum u_{(1)}cv_{(1)}\langle
u_{(2)}dv_{(2)};a,b\rangle\\
&&+\sum u_{(1)}dv_{(1)}\langle u_{(2)}cv_{(2)};a,b\rangle+\sum
u_{(1)}cdv_{(1)}\langle
u_{(2)}v_{(2)};a,b\rangle\\
&&-ucdvab+udcvab-\sum
u_{(1)}\langle u_{(2)};c,d\rangle vab\\
&\equiv&-udcvba+\sum u_{(1)}\langle u_{(2)};c,d\rangle vba+\sum
u_{(1)}v_{(1)}\langle u_{(2)}cdv_{(2)};a,b\rangle\\
&&+\sum u_{(1)}cv_{(1)}\langle u_{(2)}dv_{(2)};a,b\rangle+\sum
u_{(1)}dv_{(1)}\langle u_{(2)}cv_{(2)};a,b\rangle\\
&&+\sum u_{(1)}cdv_{(1)}\langle u_{(2)}v_{(2)};a,b\rangle
+udcvba-\sum u_{(1)}v_{(1)}\langle
u_{(2)}dcv_{(2)};a,b\rangle\\
&&-\sum u_{(1)}cv_{(1)}\langle u_{(2)}dv_{(2)};a,b\rangle-\sum
u_{(1)}dv_{(1)}\langle u_{(2)}cv_{(2)};a,b\rangle\\
&&-\sum u_{(1)}dcv_{(1)}\langle u_{(2)}v_{(2)};a,b\rangle+\sum
u_{(1)}v_{(1)}\langle u_{(2)}\langle u_{(3)};c,d\rangle
v_{(2)};a,b\rangle\\
&&+\sum u_{(1)}\langle u_{(2)};c,d\rangle
v_{(1)}\langle u_{(3)}v_{(2)};a,b\rangle-\sum u_{(1)}\langle u_{(2)};c,d\rangle vba\\
 &\equiv&\sum u_{(1)}v_{(1)}\langle
u_{(2)}[c,d]v_{(2)};a,b\rangle+\sum u_{(1)}[c,d]v_{(1)}\langle
u_{(2)}v_{(2)};a,b\rangle\\
&&+\sum u_{(1)}v_{(1)}\langle u_{(2)}\langle u_{(3)};c,d\rangle
v_{(2)};a,b\rangle+\sum u_{(1)}\langle u_{(2)};c,d\rangle
v_{(1)}\langle u_{(3)}v_{(2)};a,b\rangle\\
 &\equiv&\sum(u_{(1)}[c,d]+u_{(1)}\langle
u_{(2)};c,d\rangle)v_{(1)}\langle u_{(3)}v_{(2)};a,b\rangle\\
&&+\sum u_{(1)}v_{(1)}\langle u_{(2)}[c,d]v_{(2)};a,b\rangle+\sum
u_{(1)}v_{(1)}\langle u_{(2)}\langle u_{(3)};c,d\rangle
v_{(2)};a,b\rangle\\
&\equiv&0
\end{eqnarray*}
since $\langle x[a,b]y; c,e\rangle +\sum\langle x_{(1)}\langle
x_{(2)};a,b\rangle y; c,e\rangle=0$.

Hence, $I$ is a Gr\"obner-Shirshov basis in $mod\langle
X\rangle_{k\langle X\rangle}$.\ \ $\square$

\noindent{\bf Remark}: From the above proof, we know that for $
\widetilde{S}(V)=mod\langle X|I\rangle_{k\langle X\rangle} $, the
minimal Gr\"obner-Shirshov basis is\\
$G=\{xab-xba+\sum x_{(1)}\langle x_{(2)};a,b\rangle|x=a_{i_1}\cdots
a_{i_n}(i_1\leq \cdots\leq i_n,\ n\geq0),a>b, \
 a,b\in X\}$.

\ \

Now, by Lemma \ref{l3.5} and Theorem \ref{t7.1}, we can easily get
the following theorem.
\begin{theorem}\em(\cite{pe}, Poincare-Birkhoff-Witt) Let $\{a_i|i\in \Lambda\}$ be a totally ordered basis of
$V$. Then $\{a_{i_1}\cdots a_{i_n}|i_1\leq i_2\leq\cdots \leq i_n,\
n\geq0\}$ is a basis of $\widetilde{S}(V)$.
\end{theorem}

\noindent{\bf Acknowledgement}: The authors would like to express
their deepest gratitude to Professor L. A. Bokut for his kind
guidance, useful discussions and enthusiastic encouragement. We also
thank Professor J. M. Perez-Izquierdo for some remarks.

\ \

\end{document}